\documentclass{article}

\usepackage{graphicx}
\graphicspath{{QuiFigures/}}
\DeclareGraphicsExtensions{.png,.pdf}

\usepackage{amsmath,amsfonts}

\newcommand{\WENO}  [1][]{\ensuremath{\mathsf{WENO#1}}}

\newcommand{\CWENO} [1][]{\ensuremath{\mathsf{CWENO#1}}}

\newcommand{\DIRK}[1][]{\ensuremath{\mathsf{DIRK#1}}}
\newcommand{\IE}{\ensuremath{\mathsf{IE}}}
\newcommand{\DP}{\ensuremath{\mathsf{D3P1}}}
\newcommand{\QP}{\ensuremath{\mathsf{Q3P1}}}

\usepackage{tikz,pgfplots}
\usepackage{multirow}
\usepackage{pgfplotstable,booktabs}
\usepackage[caption=false]{subfig}

\usepackage{color}
\definecolor{Myorange}{cmyk}{0.3,0.8,1,0}
\definecolor{Myorange2}{cmyk}{0.3,0.42,1,0}
\definecolor{Mygreen}{rgb}{0.4,0.9,0.5}
\definecolor{Mygreen2}{RGB}{0,100,0}
\definecolor{Myviolet}{RGB}{148,0,211}
\definecolor{Mycyan}{cmyk}{1,0.1,0,0}
\definecolor{Myyellow}{cmyk}{0.2,0.5,1.0,0.1}
\definecolor{Myred}{rgb}{0.8,0.0,0.2}
\definecolor{Myblue}{RGB}{0,0,205}
\definecolor{darkgreen}{rgb}{0.1,0.75,0.31}

\usepackage{nicefrac}

\newcommand{\Ogrande}{\mathcal{O}}

\newcommand{\Rec}{\mathcal{R}}
\newcommand{\NumFlux}{\mathcal{F}}
\newcommand{\dx}{\mathrm{d}x}
\renewcommand{\d}{\mathrm{d}}
\def\TV{\mbox{TV}}

\newcommand{\Prec}{P_{\text{\sf rec}}}

\newcommand{\ca}[1]{\overline{#1}}
\newcommand{\DT}{\mathrm{\Delta} t}

\newcommand{\D}{\mathrm{d}}

\newcommand{\TotDer}[2]{\frac{\mathrm{d}{#1}}{\mathrm{d}{#2}}}

\title{Quinpi: integrating conservation laws with CWENO implicit methods}

\date{\today}
\author{
	Gabriella Puppo\footnote{Dipartimento di Matematica -- Sapienza, Universit\`{a} di Roma; P.le Aldo Moro, 5 -- 00185 Roma (Italy); gabriella.puppo@uniroma1.it}
	\quad Matteo Semplice\footnote{Dipartimento di Scienza e Alta Tecnologia -- Universit\`{a} dell'Insubria; Via Valleggio, 11 -- 22100 Como (Italy); matteo.semplice@uninsubria.it}
	\quad Giuseppe Visconti\footnote{Dipartimento di Matematica -- Sapienza, Universit\`{a} di Roma; P.le Aldo Moro, 5 -- 00185 Roma (Italy); giuseppe.visconti@uniroma1.it}
}

\begin{document}

\maketitle

%
%
%

\begin{abstract}
Many interesting applications of  hyperbolic systems of equations are stiff, and require the time step to satisfy restrictive stability conditions. One way to avoid small time steps is to use implicit time integration. Implicit integration is quite straightforward for first order schemes. High order schemes instead need also to control spurious oscillations, which requires limiting in space and time also in the implicit case. We propose a framework to simplify considerably the application of high order non oscillatory schemes through the introduction of a low order implicit predictor, which is used both to set up the nonlinear weights of a standard high order space reconstruction, and to achieve limiting in time. 

In this preliminary work, we concentrate on the case of a third order scheme, based on \DIRK\ integration in time and \CWENO\ reconstruction in space. The numerical tests involve linear and nonlinear scalar conservation laws.
\end{abstract}


\section{Introduction} \label{sec:intro}
Hyperbolic systems of conservation laws in one dimension can be written in the form
\begin{equation}\label{e:HypSys}
u_t+f(u)_x=0,
\end{equation}
where $u(x,t)$ is the unknown solution, and $f(u)$ is the flux function. The system is hyperbolic provided the Jacobian of $f$, $J(f)$, has real eigenvalues and a complete set of eigenvectors.

These systems model propagation phenomena, where initial and boundary data travel and interact along the eigenvectors of the system, with finite speed. The eigenvalues of $J(f)$, which depend on the solution itself when the flux is a nonlinear function of $u$, are the propagation speeds of the system. The size of the eigenvalues can span different orders of magnitude in many applications, and this fact introduces difficulties in the numerical integration of hyperbolic systems of PDE's.

In this paper, we will consider the method of lines (MOL), which is a very popular approach in the integration of hyperbolic systems. See for instance the textbook \cite{LeVeque:book}, and the classic review \cite{Shu97}, but other approaches are also effective, as in \cite{Dumbser:2008:PnPm}. In our case, we introduce a grid in the computational domain of points $x_j$. Here for simplicity we will consider a uniform grid, so all points are separated by a distance $h=x_{j}-x_{j-1}$. The computational domain $\Omega$ is thus covered with a mesh of cells $\Omega_j=[x_j-\nicefrac{1}{2}\,h,x_j+\nicefrac{1}{2}\,h]$, such that $\cup_j \Omega_j=\Omega$. 
Introducing the cell averages of the exact solution $\ca{u}(t)_j=\tfrac1h \int_{\Omega_j}u(x,t)\, \dx$, the hyperbolic system \eqref{e:HypSys} can be written as
\begin{equation} \label{e:MOL}
\TotDer{\ca{u}_j}{t} = -\frac{1}{h}\left[ f(u(x_j+\tfrac12 h,t) - f(u(x_j-\tfrac12 h,t)\right],
\end{equation}
which gives the exact evolution of the cell averages in terms of the difference of the fluxes at the cells interfaces. To transform this relation in a numerical scheme, one introduces reconstruction algorithms, $\Rec$, whose task is to estimate the value of the solution at the cell interfaces from cell averages, and numerical fluxes $\NumFlux$, which estimate the flux at the cell interfaces. Thus, the structure of a numerical scheme in the method of lines approach can be written as
\begin{enumerate}
\item Define a reconstruction algorithm $\Rec(\{\ca{u}\}_j)$ such that the exact solution $u$ at time $t$ is approximated with $u(x,t) = \sum_j R_j(x; t)\chi_j(x)$, where $\chi_j$ is the characteristic function of the interval $\Omega_j$, and $R_j$ is the restriction of $\Rec$ to the interval $\Omega_j$. Typically, $R_j$ is a polynomial of degree $d_j$, which changes in time because it is defined starting from the time-dependent cell averages. In this paper, we will assume that $d_j\equiv d$ is constant.
\item Compute the boundary extrapolated data (BED) at the cell interfaces $u^+_{j+\nicefrac12}=R_{j+1}(x_j+\tfrac{h}{2})$ and $u^-_{j+\nicefrac12}=R_{j}(x_j+\tfrac{h}{2})$. Recall that the two BED's computed at the same interface are different, with $u^+_{j+\nicefrac12}-u^-_{j+\nicefrac12}=\Ogrande(h^{d+1})$ when the flow is smooth enough.
\item Choose a smooth enough numerical flux function $\NumFlux(a,b)$ such that $\NumFlux(a,a)=f(a)$, with the stability property that $\NumFlux$ be an increasing function of the first argument and a decreasing function of the second argument.
\item Then the solution of the PDE is approximated by the solution of the system of ODE's
\begin{equation} \label{e:spaceApprox}
\TotDer{\ca{u}_j}{t} = - \frac{1}{h}\left[ \NumFlux_{j+\nicefrac12}-\NumFlux_{j-\nicefrac12}\right],
\end{equation}
with $\NumFlux_{j+\nicefrac12}=\NumFlux(u^-_{j+\nicefrac12},u^+_{j+\nicefrac12})$.
\end{enumerate}

It is well known that for explicit schemes the time step $\DT$ must satisfy the CFL condition, namely $\lambda=\DT/h\leq 1/\max_u|f'(u)|$. Thus the numerical speed $1/\lambda$ must be {\em faster} than all the waves present in the system. Since the error in the numerical solution depends on the difference between the numerical and the actual speed, it follows that in explicit schemes the waves that are better approximated  are the fastest waves.

Several systems of hyperbolic conservation laws are characterized by waves with very different speeds, and in many applications the phenomenon of interest travels with slow speeds, while the fastest waves which impose the CFL condition do not need to be accurately represented.
Low Mach problems in gas dynamics and kinetic problems close to equilibrium provide good examples.

Low Mach problems arise in gas dynamics when the flow is close to incompressibility. In these cases the actual speed of the gas is much slower than the acoustic waves, and if one is interested in the movement of the gas, accuracy in the propagation of sound is irrelevant. A huge literature has developed for these problems, here we mention just a few pioneering works, \cite{2010DellacherieLowMach,2011DegondTang} and some more recent developments \cite{2017AbbateAllSpeed,2017DimarcoLoubere_LowMachAP,2018BoscarinoScandurra,2017Tavelli_SemiImplicitAllMach}.

A second example in which fast velocities constraint the CFL condition, but the signals carried by them do not need to be accurately represented, occurs in kinetic problems, especially close to equilibrium. In kinetic problems the typical evolution equation has the form
\begin{equation}
f_t + v \cdot \nabla_x f= Q(f,f), \label{e:kinetic}
\end{equation}
where $f=f(x,t,v)$ is the distribution function for particles located at $x$, with velocity $v$ at time $t$, and $Q(f,f)$ is the collision term, accounting for the interaction between particles during which the microscopic speeds $v$ are modified. When the flow is close to equilibrium, the relevant phenomena travel with macroscopic speeds which have magnitude of order $\int v \,f(x,t,v) \D v$ and $[\int v^2 \,f(x,t,v) \D v]^{\nicefrac12}$. Still, the CFL condition is based on the fastest microscopic speeds, which are typically much larger. Several attempts have been proposed to go around this restriction, see the review \cite{DimarcoPareschiActa}, and \cite{MicroMacro}. Note that the microscopic speeds appear in the convective term only, and convection is {\em linear} in kinetic problems. Thus it is feasible to integrate these equations with implicit schemes, taking advantage of the linearity of the convective terms, as in \cite{MiMe}.

Another setting in which convection is linear, and therefore more amenable to implicit integration, arises with relaxation systems of the form proposed in \cite{1995JinXin_Relaxation}. Relaxation leads to systems of PDE's  which have a kinetic form with linear transport, as in \eqref{e:kinetic}. Besides kinetic problems, diffusive relaxation leads to fast relaxation microscopic speeds, which again prompts the need for implicit time integration.

In this paper, we will propose numerical schemes for the implicit integration in time of hyperbolic systems of equations. This is a preliminary work. For the time being, we will focus on third order implicit schemes for scalar conservation laws. As noted in \cite{2020Arbogast}, there are three levels of non-linearity in high order implicit methods for conservation laws, which make implicit integration particularly challenging. The first level consists in the non-linearity of the flux function, which is due to the physical structure of the model, and therefore is unavoidable. The other sources of non-linearity are due to the need to prevent spurious oscillations, arising with high order numerical schemes. Even in the explicit case, non oscillatory high order schemes must use nonlinear reconstructions. As is well known,  nonlinear reconstructions in space are needed even in implicit schemes. Moreover, high order time integrators are typically based on a polynomial approximation of the time derivative. Thus, a nonlinear limiting is needed also in the time reconstruction of the derivative. These problems have been addressed by several authors. We mention work on second order schemes in \cite{2003DurasaisamyBaeder,2007DurasaisamyBaeder}, where limiters are applied in space and time simultaneously and TVD estimates are derived. An interesting discussion on TVD bounds in space and time can be found in \cite{2001Forth}, and in the classic paper \cite{1984Harten}. See also \cite{2001Gottlieb}. A fully nonlinear third order implicit scheme, which is limited in both space and time simultaneously can be found in \cite{2020Arbogast}. 

The next section, \S~\ref{s:motivation}, contains a discussion of revisited TVD bounds which are at the basis of our approach, and a summary of the choices leading to the final scheme. Section \ref{sec:scheme} is the main part of the work, and it offers a complete discussion of the bricks composing the proposed Quinpi algorithm. Next, section \ref{sec:numerics} documents the properties of the scheme with a selection of numerical tests for scalar equations. We end with conclusions and a plan for future work in \S~\ref{sect:conclusion}.

\section{Motivation}\label{s:motivation}

To study TVD conditions for implicit schemes, we consider the linear advection equation $u_t+au_x=0$, with $a>0$ and the upwind scheme. In this section, $u_j^n$ is the numerical approximation to the exact solution, at one of the space-time grid points, namely $u_j^n\simeq u(x_j,t^n)$, with $x_j\in \Omega$, is a grid point in the space mesh, and $t^n=n\Delta t$, $\Delta t$ being the time step. The mesh ratio is $\lambda=\Delta t/h$.

In the explicit case, we have
\[
u_j^{n+1}=u_j^n -\lambda a(u_j^n-u_{j-1}^n)
\]
and the total variation of $u^{n+1}$ is
\begin{align*}
\TV(u^{n+1})& = \sum_j\left| (1-\lambda a)(u^n_j-u^n_{j-1})
      +\lambda a(u_{j-1}^n-u_{j-2}^n)\right| \\
      & \leq |1-\lambda a| \sum_j\left|u^n_j-u^n_{j-1}\right|
      +|\lambda a|\sum_j\left|u^n_{j-1}-u^n_{j-2}\right|.
\end{align*}
Then, using the fact that $a>0$, periodic boundary conditions, or compact support of the solution, and the CFL condition, one has that $\TV(u^{n+1})\leq \TV(u^n)$.

For the implicit upwind scheme instead, we start from
\[
u_j^n= u_j^{n+1}+\lambda a(u_j^{n+1}-u_{j-1}^{n+1}).
\]
Now
\begin{align*}
\TV(u^{n})& = \sum_j\left| (1+\lambda a)(u^{n+1}_j-u^{n+1}_{j-1})
      -\lambda a(u_{j-1}^{n+1}-u_{j-2}^{n+1})\right| \\
      & \geq |1+\lambda a| \sum_j\left|u^{n+1}_j-u^{n+1}_{j-1}\right|
      -|\lambda a|\sum_j\left|u^{n+1}_{j-1}-u^{n+1}_{j-2}\right|.
\end{align*}
Since $a>0$, and applying again periodic boundary conditions, we find $\TV(u^{n+1})\leq \TV(u^n), \ \forall \lambda$. Thus, the implicit upwind scheme is not only unconditionally stable, but it is also unconditionally Total Variation non increasing.

We now consider the explicit Euler scheme with second order space differencing, using a piecewise limited reconstruction. The numerical solution now is
\[
u_j^{n+1}=u_j^n -\lambda a\left(1-\tfrac12\Phi(\theta_{j-1})\right)\left(u_j^n-u_{j-1}^n\right)-\tfrac12 \lambda a\Phi(\theta_j) \left(u_{j+1}^n-u_{j}^n\right),
\]
where $\Phi(\theta)$ is the limiter, and the quantity $\theta_j=\tfrac{u_j-u_{j-1}}{u_{j+1}-u_{j}}$, so that the limited slope is $\sigma_j=\Phi(\theta_j)(u_{j+1}-u_{j})$.

To prove under what conditions the scheme is TVD, following Harten \cite{1983Harten_HighResolution},  and Sweby \cite{1984Sweby_TVD}, one rewrites this scheme as
\[
u_j^{n+1}=u_j^n -C_{j-1}\left(u_j^n-u_{j-1}^n\right),
\]
with
\[
C_{j-1}=\lambda a\left[1-\tfrac12 \left( \Phi(\theta_{j-1})-\frac{\Phi(\theta_{j})}{\theta_{j}}\right)\right].
\]
Computing the Total Variation of the scheme above, one finds that the scheme is TVD provided
\[
C_j \geq 0, \qquad \mbox{and}\quad 1-C_j\geq 0.
\]
Since $a>0$, the first condition holds provided
\[
\alpha_j = 1-\tfrac12 \left( \Phi(\theta_{j-1})-\frac{\Phi(\theta_{j})}{\theta_{j}}\right) \geq 0.
\]
From this, one recovers the familiar restrictions in the choice of the limiter function, namely $0\leq \tfrac{\Phi(\theta)}{\theta}\leq 2$ and $0\leq \Phi(\theta)\leq 2$, see also LeVeque \cite{LeVeque:book}. Since $\alpha_j\geq 0$, the second condition on $C_{j}$ is satisfied provided CFL holds.

How do these estimates change in the implicit case? Now we have
\[
u_j^n= u_j^{n+1}+\lambda a\left(1-\tfrac12\Phi(\theta_{j-1})\right)\left(u_j^{n+1}-u_{j-1}^{n+1}\right)+\tfrac12 \lambda a\Phi(\theta_j) \left(u_{j+1}^{n+1}-u_{j}^{n+1}\right),
\]
and introducing the same quantity $C_j$ we saw before, this expression rewrites as
\[
u_j^n= u_j^{n+1}+C_{j-1}\left(u_j^{n+1}-u_{j-1}^{n+1}\right).
\]
The Total Variation of this scheme is
\begin{align*}
\TV(u^{n})& = \sum_j\left| (1+C_{j-1})(u^{n+1}_j-u^{n+1}_{j-1})
      -C_{j-2}(u_{j-1}^{n+1}-u_{j-2}^{n+1})\right| \\
      & \geq \sum_j|1+C_{j-1}| \left|u^{n+1}_j-u^{n+1}_{j-1}\right|
      -\sum_j|C_{j-2}|\left|u^{n+1}_{j-1}-u^{n+1}_{j-2}\right|,
\end{align*}
where we have used the inverse triangle inequality. Applying periodic boundary conditions, and assuming that
\[
C_j \geq 0, \qquad \mbox{and}\quad 1+C_j\geq 0,
\]
one finds that the scheme is TVD. The important point for the development in this work is that the first condition gives {\em exactly the same restrictions} on the choice of the limiter function we find in the explicit case. The second condition instead is always satisfied. This means that the implicit scheme with piecewise linear reconstruction in space is TVD for all $\lambda$, provided the limiter function satisfies the usual bounds holding for explicit schemes. 

The limiter is subject also to accuracy constraints. We achieve second order accuracy if the limited slope is a first order approximation of the exact slope. Namely, if $U(x)$ is a smooth function with cell averages $\ca{u}_j$, second order accuracy requires that
\[
U(x_{j+\nicefrac12})-\left[ \ca{u}_j+\tfrac12 \sigma_j\right]=\Ogrande(h)^2.
\] 
Substituting the limited slope and expanding around $U(x_j)=u_j$, one finds
\[
U(x_{j+\nicefrac12})-u^-_{j+\nicefrac12}=\tfrac{h}{2}(1-\Phi(1))+\Ogrande(h)^2.
\] 
As in the explicit case, we see that second order accuracy requires $\Phi(1)=1$, together with at least Lipshitz continuity for $\Phi$. This means that $\Phi$ and $\theta$ can be computed using auxiliary functions that must be at least first order approximations of the unknown function $u(x,t^{n+1})$. For these reasons, in this work we will
\begin{itemize}
\item Build implicit schemes for stiff hyperbolic systems with unconditional stability, based on diagonally implicit Runge Kutta schemes (\DIRK).
\item Prevent spurious oscillations in space using {\em exactly the same} techniques for {\em space limiting} we are familiar with. In fact, as pointed out in \cite{2001Gottlieb}, Runge Kutta schemes can be seen as combinations of Euler steps. This is true also in the implicit case. We will strive to ensure that each implicit Euler step composing the whole Runge Kutta step prevents spurious oscillations.  In other words, we will rewrite the PDE as a system of ODE's in the cell averages, with a non oscillatory right hand side.
\item The non-linearities in the space reconstruction will be tackled using first order predictors, that can be computed without limiting, because, as we have recalled, they are unconditionally TVD. They allow to compute each \DIRK\ stage by solving a system which is nonlinear only because of the nonlinear flux function, and, at the same time, they will constitute a low order non oscillatory approximation of the the solution in its own right.
\item As is known, this is not enough to prevent spurious oscillations, when high order space reconstructions are coupled with high order time integrators. A high order Runge Kutta scheme can be viewed as a polynomial reconstruction in time through natural continuous extensions \cite{NorsettWanner:1981,Zennaro}. Not surprisingly, limiting is needed also in time, as pointed out in \cite{2007DurasaisamyBaeder} and more recently \cite{2020Arbogast}. Unlike these authors, however, we limit the solution in time applying the limiter directly on the computed solution, blending the accurate solution with the low order predictor, without coupling space and time limiting.
\end{itemize}

We will name the resulting schemes Quinpi, for \CWENO\ Implicit.
In this paper, we will introduce only third order Quinpi schemes. We plan to extend the construction to include higher order schemes, BDF extensions and applications in future works.

\section{Quinpi finite volume scheme} \label{sec:scheme}

This section is devoted to the description of the implicit numerical solution of a one-dimensional scalar conservation law of the form~\eqref{e:HypSys} with the Quinpi approach. 
As pointed out in Section~\ref{sec:intro}, we focus on the finite volume framework employing the method of lines, which leads to the system of ODEs~\eqref{e:spaceApprox} for the evolution of the approximate cell averages. 

In the following subsections, we will initially describe the space reconstruction for a function $u(x)$, with limiters based on a predicted solution $p(x)$, where $p$ should be at least a first order approximation of $u$. We will concentrate on the third order \CWENO\ reconstruction \cite{LPR:00:SIAMJSciComp}. For the general \CWENO\ algorithm, see \cite{CPSV:cweno}. More details, improvements and extensions can be found in \cite{SCR:CWENOquadtree,CWENOandaluz,CSV19:cwenoz,CS:epsweno,DBSR:ADER_CWENO,Kolb:14}. Other non oscillatory reconstructions can also be used, such as \WENO, see the classic review \cite{Shu97}. and its extensions as \cite{Balsara:AOWENO}.

Next, we will consider the integration in time, with a \DIRK[3]\ scheme, describing the interweaving of predictor and time advancement of the solution. Finally, we will introduce the limiting in time, which consists in a nonlinear blending between the predicted and the high order solution, and the conservative correction that is needed at the end of the time step.

\subsection{Space approximation: the \CWENO\ reconstruction procedure} \label{sec:rec}

Central Weighted Essentially Non Oscillatory \CWENO\ schemes are a class of high-order numerical methods to reconstruct accurate
and non-oscillatory point values of a function $u$ starting from  the knowledge of its cell averages. The main characteristic of \CWENO\ reconstructions is that they are uniformly accurate within the entire cell, \cite{CPSV:cweno}. 

For the purpose of this work, here we briefly recall the definition of the \CWENO\ reconstruction procedure restricting the presentation to the one-dimensional third order scalar case on a uniform mesh. In order to reconstruct a function $u$ at some $x\in \Omega$ and at a fixed time $t$, we consider as given data the cell averages $\ca{u}_j$ of $u$ at time $t$ over the cells $\Omega_j$ of a grid, which is a uniform discretization of $\Omega$.

A third-order \CWENO\ reconstruction is characterized by the use of an optimal polynomial of degree $2$ and two linear polynomials. Let $P_{j,L}^{(1)}$ and $P_{j,R}^{(1)}$ be the linear polynomials
\begin{equation} \label{e:plinear}
\begin{aligned}
P_{j,L}^{(1)} &= \overline{u}_j + \frac{\overline{u}_j-\overline{u}_{j-1}}{h} (x-x_j) \\
P_{j,R}^{(1)} &= \overline{u}_j + \frac{\overline{u}_{j+1}-\overline{u}_{j}}{h} (x-x_j),
\end{aligned}
\end{equation}
and let us write the optimal polynomial of degree $2$ as
\begin{equation} \label{e:popt}
\begin{aligned}
P_j^{(2)}(x) &= a + b(x-x_j) + c(x-x_j)^2 \\
a &= \frac{-\overline{u}_{j+1} +26\overline{u}_j -\overline{u}_{j-1}}{24} \\
b &= \frac{\overline{u}_{j+1} -\overline{u}_{j-1}}{2h} \\
c &= \frac{\overline{u}_{j+1} -2\overline{u}_j +\overline{u}_{j-1}}{2h^2}.
\end{aligned}
\end{equation}
All polynomials $P_{j,L}^{(1)}$, $P_{j,R}^{(1)}$ and $P_j^{(2)}$ interpolate the data in the sense of cell averages. 
We also introduce the second degree polynomial $P_{j,0}$ defined as
\begin{equation} \label{e:p0}
\begin{aligned}
P_{j,0}(x) &= \frac{1}{C_0} \left( P_j^{(2)}(x) - C_L P_{j,L}^{(1)} - C_R P_{j,R}^{(1)} \right) = A + B(x-x_j) + C(x-x_j^2) \\
A &= \frac{a}{C_0} - \frac{C_L+C_R}{C_0} \overline{u}_j \\
B &= \frac{(1-2C_R)\overline{u}_{j+1} -(2C_L-2C_R)\overline{u}_j +(2C_L-1)\overline{u}_{j-1}}{2C_0h} \\
C &= \frac{c}{C_0}, 
\end{aligned}
\end{equation}
where $C_0, C_L, C_R \in (0,1)$ with $C_0+C_L+C_R = 1$ are the linear or optimal coefficients. Note that this polynomial reproduces the underlying data with first order accuracy.

The Jiang-Shu smoothness indicators \cite{JiangShu:96} of the polynomials related to the cell $\Omega_j$ are defined by
\begin{equation} \label{e:ind}
	I[P^{(r)}] := \sum_{i=1}^r  h^{2i-1} \int_{\Omega_j} \left(\frac{\d^i}{\dx^i} P^{(r)}(x)\right)^2 \dx.
\end{equation}
Then, in our case, they reduce to
\begin{equation} \label{e:indicators}
\begin{aligned}
I_{j,L} = (\overline{u}_j -\overline{u}_{j-1})^2, \quad I_{j,R} = (\overline{u}_{j+1} -\overline{u}_{j})^2, \quad I_{j,0} = b^2 h^2 + \frac{52}{3} c^2 h^4.
\end{aligned}
\end{equation}
From these, the nonlinear weights are defined as
\begin{equation} \label{e:nonlinw}
\tilde{\omega}_{j,k} = \frac{C_{k}}{(\epsilon_x+I_{j,k})^\tau}, \quad \omega_{j,k} = \frac{\tilde{\omega}_{j,k}}{\tilde{\omega}_{j,0}+\tilde{\omega}_{j,L}+\tilde{\omega}_{j,R}}, \quad k = 0,L,R.
\end{equation}
Following \cite{CS:epsweno,CSV19:cwenoz,Kolb:14}, we will always set $\epsilon_x=h^2$ and $\tau=2$.
The reconstruction polynomial $P_{j,\text{rec}}$ is
\begin{equation} \label{e:prec}
P_{j,\text{rec}}(x) = \omega_{j,0} P_{j,0}(x) + \omega_{j,L} P_{j,L}^{(1)}(x) + \omega_{j,R} P_{j,R}^{(1)}(x).
\end{equation}
Note that if $\omega_{j,k}=C_k$, $k=0,L,R$, then the reconstructed polynomial $P_{j,\text{rec}}$ coincides with the optimal polynomial $P_j^{(2)}$, and the reconstruction is third order accurate in the whole cell.

If the nonlinear weights satisfy
\begin{equation} \label{e:weightsDifference}
	C_k-\omega_{j,k}=\Ogrande(h),
\end{equation}
then \CWENO\ boosts the accuracy of the reconstruction polynomial $\Prec$ to the highest possible accuracy, $3$ in this case. This occurs when the stencil containing the data is smooth. 

Once the reconstruction is known, we can compute the boundary extrapolated data as
\begin{equation}\label{e:BED}
u^-_{j+\nicefrac12}= P_{j,\text{rec}}(x_{j+\nicefrac12}) 
\qquad 
u^+_{j+\nicefrac12}= P_{j+1,\text{rec}}(x_{j+\nicefrac12}) 
\end{equation}
In Quinpi methods, we suppose we are given a predictor $p(x)$ alongside the solution $u(x)$, where the predictor is at least a first order approximation of the solution. We compute the smoothness indicators \eqref{e:ind} using the predictor. Thus, the weights $\omega_{j,k}$ depend only on the predictor. In this fashion, the boundary extrapolated data \eqref{e:BED} are {\em linear functions} of the solution $u$ with coefficients that are not constant across the grid.

\subsection{Time approximation: diagonally implicit Runge-Kutta} \label{sec:time}

The space reconstruction algorithm allows to compute boundary extrapolated data of the solution $u$ of the conservation law~\eqref{e:HypSys} at each cell interface. 
Next we pick a consistent and monotone numerical flux function $\mathcal{F}$. For example, in this work, we will apply the Lax Friedrichs numerical flux
\begin{equation} \label{e:lxf}
\mathcal{F}(u^-,u^+) = \frac12 \left( f(u^{+}) + f(u^{-}) - \alpha (u^{+} - u^{-}) \right), 
\end{equation}
with $\alpha = \max_u|f'(u)|$.

This completely defines the system of ODE's \eqref{e:spaceApprox}.
This system  needs to be approximated in time by means of a time integration scheme. Here, we focus on Diagonally Implicit Runge-Kutta (\DIRK) methods, with $s$ stages and general Butcher tableau
\begin{equation}
\label{e:tableau}
\begin{array}{c|cccc}
c_1 & a_{11} & 0 & \dots & 0 \\[1.5ex]
c_2 & a_{21} & a_{22} & \dots & 0 \\[1.5ex]
\vdots & \vdots & \vdots & \ddots & \vdots \\[1.5ex]
c_s & a_{s1} & a_{s2} & \dots & a_{ss} \\[1ex]
\hline
&&&&\\[-3ex]
& b_1 & b_2 & \dots & b_s
\end{array}
\end{equation}
having the property $a_{k,\ell}=0$, for each $k<\ell$. 

Discretization of~\eqref{e:spaceApprox} with a \DIRK\ method leads to the fully discrete scheme
\begin{equation} \label{e:scheme}
	\overline{u}_j^{n+1} = \overline{u}_j^{n} - \frac{\Delta t}{h} \sum_{k=1}^s b_k \left[ \mathcal{F}_{j+\frac12}^{(k)} - \mathcal{F}_{j-\frac12}^{(k)} \right], \quad n\geq 0, \ j=1,\dots,N
\end{equation}
where we recall that $h$ is the mesh spacing, $\Delta t$ is the time step and $\overline{u}_j^n \approx \ca{u}_j(t^n)$. Finally, $\mathcal{F}_{j+\frac12}^{(k)}=\mathcal{F}(u_{j+\frac12}^{-,(k)},u_{j-\frac12}^{+,(k)})$, where the boundary extrapolated data $u_{j+\frac12}^{-,(k)}$, $u_{j-\frac12}^{+,(k)}$, are reconstructions at the cell boundaries of the stage values
\begin{equation} \label{e:stagevalue}
	\overline{u}_j^{(k)} = \overline{u}_j^n - \frac{\Delta t}{h} \sum_{\ell=1}^{k} a_{k,\ell} \left[ \mathcal{F}_{j+\frac12}^{(\ell)} - \mathcal{F}_{j-\frac12}^{(\ell)} \right], \quad k=1,\dots,s,
\end{equation}
which are approximations at times $t^{(k)} = t^n + c_k\Delta t$.
We point out that $\mathcal{F}_{j+\frac12}^{(k)}$ depends on the $k-$th stage value, and thus the computation of each stage is implicit but independent from the following ones.

The \DIRK\ scheme used in this work is
\begin{equation}
\label{e:dirk3}
\begin{array}{c|ccc}
\lambda & \lambda & 0 & 0 \\[1.5ex]
\frac{(1+\lambda)}{2} & \frac{(1-\lambda)}{2} & \lambda & 0 \\[1.5ex]
1 & -\frac32 \lambda^2 + 4\lambda - \frac14 & \frac32\lambda^2-5\lambda+\frac54 & \lambda \\[1ex]
\hline
&&&\\[-3ex]
& -\frac32 \lambda^2 + 4\lambda - \frac14 & \frac32\lambda^2-5\lambda+\frac54 & \lambda
\end{array}
\end{equation}
with $\lambda=0.4358665215$, see
\cite{Alexander1977}.

The fully discrete set of equations~\eqref{e:scheme}-\eqref{e:stagevalue} has two sources of non-linearity when solved with high-order schemes: one arises from the physics of the system~\eqref{e:HypSys} when the flux function $f$ is nonlinear, and cannot be avoided; the other one arises from the computation of the boundary extrapolated data using a high order \CWENO\ or \WENO\ reconstruction, which is nonlinear even for linear problems.  Therefore, even for a linear PDE, the resolution of~\eqref{e:scheme}-\eqref{e:stagevalue} requires a nonlinear solver. 
Typically, one uses Newton's algorithm, 
which requires the computation of the Jacobian of the scheme, which depends on the nonlinear weights~\eqref{e:nonlinw} and on the oscillation indicators~\eqref{e:ind}-\eqref{e:indicators}, resulting in a prohibitive computational cost.

In the next subsection, we propose a way to circumvent the non-linearity with a high-order reconstruction procedure relying on a predictor.

\subsection{Third-order Quinpi approach}\label{sec:fluximplicit}

A prototype of an implicit scheme for \WENO\ reconstructions based on a predictor was developed by Gottlieb, Mullen and Ruuth in~\cite{Gottlieb:iWENO:2006}. The method relies on the idea of a predictor-corrector approach to avoid the non-linearity of the reconstruction. The solution of an explicit scheme is used as predictor in order to compute the nonlinear weights of \WENO\ within the high-order implicit scheme, which is used as corrector.

In~\cite{Gottlieb:iWENO:2006} the approximation with the explicit predictor is computed in a single time step, namely without performing several steps within the Courant number. On the contrary, in the Quinpi approach, the non-linearity arising from the nonlinear weights of \CWENO\ is circumvented by computing an approximation of the solution at each intermediate time $t^{(k)} = t^n + c_k \Delta t$ with an implicit, but linear, low-order scheme, with which the nonlinear weights of \CWENO\ are predicted at each stage.

Once the weights are known, a correction of order $3$ is obtained by employing a \DIRK\ method of order three coupled with the third-order \CWENO\ space reconstruction, with the weights computed from the predictor. 

In this way, the complete scheme is linear with respect to the space reconstruction. In this context, by {\em linear} we mean that, for a linear conservation law, the solution can be advanced by a time step solving a sequence of $s$ narrow-banded linear systems. However, the scheme overall is nonlinear with respect to its initial data because the entries of the linear systems' matrices depend non-linearly on the predicted solution through~\eqref{e:nonlinw} and~\eqref{e:indicators}. When $f(u)$ is not linear, the systems become nonlinear but only through the flux function.

Clearly, an implicit predictor is more expensive to compute than the explicit predictor proposed in \cite{Gottlieb:iWENO:2006}. But using an implicit predictor has a double advantage. First, the predictor itself is stable, and this allows to have a reliable prediction of the weights, even for high Courant numbers. Second, at the end of the time step, the predictor itself is a reliable, stable and non oscillatory low order solution, with which we will blend the high order solution to obtain the time limiting required by high order time integrators.

\subsubsection{First-order predictor: composite implicit Euler}

We solve~\eqref{e:scheme}-\eqref{e:stagevalue}-\eqref{e:lxf} within the time step $\Delta t$ using an $s$ stages composite backward Euler scheme in time, where $s$ is the number of stages in the \DIRK\ scheme \eqref{e:tableau}. In other words, we apply the backward Euler scheme $s$ times in each time step. The $k$-th substep advances the solution from $t^n+\tilde{c}_{k-1} \Delta t$ to $t^n+\tilde{c}_{k} \Delta t$, with $\tilde{c}_0=0$, where the coefficients $\tilde{c}_k$, $k=1,\dots,s$, are the ordered abscissae of the \DIRK~\eqref{e:tableau}. Overall, this is equivalent to implementing a \DIRK\ scheme with Butcher tableau given by
\begin{equation}
\label{e:tableauIE}
\begin{array}{c|cccc}
\tilde{c}_1 & \tilde{c}_1 & 0 & \dots & 0 \\[1.5ex]
\tilde{c}_2 & \tilde{c}_1 & \tilde{c}_2-\tilde{c}_1 & \dots & 0 \\[1.5ex]
\vdots & \vdots & \vdots & \ddots & \vdots \\[1.5ex]
\tilde{c}_s & \tilde{c}_1 & \tilde{c}_2-\tilde{c}_1 & \dots & \tilde{c}_s-\tilde{c}_{s-1} \\[1ex]
\hline
&&&&\\[-3ex]
& \tilde{c}_1 & \tilde{c}_2-\tilde{c}_1 & \dots & \tilde{c}_s-\tilde{c}_{s-1}
\end{array}
\end{equation}
In space, we use piecewise constant reconstructions from the cell averages. At the final stage, we thus obtain a first order stable non oscillatory approximation of the solution at the time step $t^{n+1}$, which we will call $u^{\IE, n+1}$. A similar low order composite backward Euler scheme is also employed in the third order scheme of Arbogast et al.~\cite{2020Arbogast}.

The resulting scheme also provides first-order approximations $\ca{u}_j^{\IE,(k)}$ of the solution at the intermediate times $t^{(k)}=t^n + \tilde{c}_k \Delta t$, for $k=1,\dots,s$.
At each stage therefore one needs to solve the
nonlinear system 
\begin{equation} \label{e:nonlinearIE}
\begin{aligned}
	G_j(\overline{U}^{\IE,(k)}) :=& \frac{\theta_k \Delta t}{h} \left[ \frac12 f(\overline{u}_{j+1}^{\IE,(k)}) - \frac12 f(\overline{u}_{j-1}^{\IE,(k)}) - \frac{\alpha}{2} \overline{u}_{j+1}^{\IE,(k)} + \right. \\
	&\left. + \left( \alpha + \frac{h}{\theta_k \Delta t} \right) \overline{u}_j^{\IE,(k)} - \frac{\alpha}{2} \overline{u}_{j+1}^{\IE,(k)} \right] - \overline{u}_j^{\IE,(k-1)} = 0,
\end{aligned}
\end{equation}
where, $\theta_k:=\tilde{c}_{k}-\tilde{c}_{k-1}$, $\overline{U}^{\IE,(k)}:=\{\ca{u}_j^{\IE,(k)}\}_j$ and $\overline{u}_j^{\IE,(0)}:=\overline{u}_j^{n}$. We use Newton's method. Note that the system is nonlinear only through the flux function.

\subsubsection{\CWENO\ third-order correction}

Once the low order predictions $\ca{u}_j^{\IE,(k)}$ are known at all times $t^{(k)} = t^n+\tilde{c}_k\Delta t$, we correct the accuracy of the solution by solving~\eqref{e:scheme}-\eqref{e:stagevalue}-\eqref{e:lxf} using the third-order \DIRK\ \eqref{e:dirk3} in time and the third-order \CWENO\ reconstruction in space, with the weights $\omega_j^{(k)}$ at the $k$-th stage computed through the predictor $\ca{u}_j^{\IE,(k)}$. Thus the boundary extrapolated data will be given by
\begin{equation} \label{e:reconstructionsFI}
\begin{aligned}
u_{j+\frac12}^{-,(k)} = P_{j,\text{rec}}(x_{j+\frac12}) &= \sum_{\ell=-1}^1 W_{j,\ell}^{-} \overline{u}_{j+\ell}^{(k)} \\
u_{j+\frac12}^{+,(k)} = P_{j+1,\text{rec}}(x_{j+\frac12}) &= \sum_{\ell=-1}^1 W_{j,\ell}^{+} \overline{u}_{j+1+\ell}^{(k)},
\end{aligned}
\end{equation}
where the weights $W_j^{\pm}$ depend only on $\ca{u}_j^{\IE,(k)}$ and are constant with respect to $\ca{u}^k$.
Finally, from~\eqref{e:scheme}~\eqref{e:stagevalue}~\eqref{e:lxf}, at each stage, we solve the system
\begin{align*}
G_j(\overline{U}^{(k)}):=& \overline{u}_j^{(k)}- \overline{u}_j^n \\
& + \frac{a_{k,k}\Delta t}{h} \left[ \mathcal{F}_{j+1/2}^{(k)}-\mathcal{F}_{j-1/2}^{(k)} \right] + \frac{\Delta t}{h} \sum_{\ell=1}^{k-1} a_{k,\ell} [\mathcal{F}_{j+1/2}^{(\ell)}-\mathcal{F}_{j-1/2}^{(\ell)}] = 0.
\end{align*}
We observe that the non-linearity of $G_j$ is only due to the flux function, and not to the \CWENO\ reconstruction, which uses the predictor $\ca{u}_j^{\IE,(k)}$ to compute the nonlinear weights. 

In the following, the numerical solution obtained applying the \DIRK3 scheme, with the third order \CWENO\ reconstruction exploiting the $u^\IE$ predictor will be called \DP\ scheme.

\subsubsection{Nonlinear blending in time}

The solution $u^{\DP}$ obtained with the \DP\ scheme is third order accurate, and has control over spurious oscillations, thanks to the limited space reconstruction described above. However, this solution may still exhibit oscillations, because it is not limited in time. We discuss in this section the definition of the time limited solution.


It is easy to associate a Continuous Extension (CE) to a Runge-Kutta scheme. For example, when the abscissae $c_k$, $k=1,\dots,s$, are distinct, following~\cite{NorsettWanner:1981} one can construct a polynomial $P(t)$ such that $P(t^n)=u^n$ and $P^\prime(t^n+c_k\Delta t)=K_k$, where the $K_k$'s are the RK fluxes of the Runge-Kutta scheme and $u^n$ is the solution at time $t^n$. The polynomial $P(t)$ is such that $P(t^n+\Delta t)=u^{n+1}$ and it provides a way to interpolate the numerical solution at any point $t^n+\gamma\Delta t$ for $\gamma\in(0,1)$. If, for a particular Runge-Kutta scheme, a Natural Continuous Extension in the sense of~\cite{Zennaro} exists, one could use that, but for this work we do not require the extensions to be natural, since we use them only as a device to assess the smoothness of the solution in the time step.

In our case, the CE of the \DIRK3 scheme we are using defines a polynomial extension of degree 3 in the sense of~\cite{NorsettWanner:1981}, which will be called $P^{(3)}_t$. Instead, we name $P^{(1)}_t$ the polynomial extension underlying the composite implicit Euler~\eqref{e:tableauIE}.

We define the limited in time solution $u^B(t)$ in a \CWENO\ fashion as
\begin{equation} \label{e:TimeBlended_general}
u^{B}_j(t)=\frac{\omega_j^{H,n}}{C_H}
\left( P^{(3)}_t(t)-C_L P^{(1)}_t(t)\right)
+\omega_j^{L,n}P^{(1)}_t(t), \quad t\in[t^n,t^{n+1}],
\end{equation}
with $\omega_j^{H,n}+\omega_j^{L,n}=1$.
In this and in all subsequent equations, $H$ and $L$ stand for high and low order quantities, respectively. The coefficients $C_L$ and $C_H$ are such that $C_L, C_H\in(0,1)$, with $C_L+C_H=1$. We observe that, by a property of the CE polynomials, at time $t^{n+1}$ we have
\begin{equation}\label{e:TimeBlended}
u^{B,n+1}_j=\frac{\omega_j^{H,n}}{C_H}
\left( u_j^{\DP,n+1}-C_Lu_j^{\IE,n+1}\right)
+\omega_j^{L,n}u_j^{\IE,n+1}.
\end{equation}
Equation~\eqref{e:TimeBlended} describes a nonlinear blending between the low order solution $u^{\IE}$ and the high order solution $u^{\DP}$ at time $t^{n+1}$. We notice that if $\omega_j^{H,n}=C_H$, and consequently $\omega_j^{L,n}=C_L$, then $u^{B,n+1}_j=u_j^{\DP,n+1}$ and the blending selects the solution of the high order scheme. Instead, if $\omega_j^{L,n}=1$, $u^{B,n+1}_j=u_j^{\IE,n+1}$ and the blending selects the solution of the low order scheme.

The weights $\omega^L$ and $\omega^H$ must be thus designed in order to privilege the high order solution when it is not oscillatory, and the low order solution otherwise. In the following we discuss their definition, which is also carried out as in \CWENO\ and relies on suitable regularity indicators. In fact, we define
\begin{equation}\label{e:TimeWeights}
\omega_j^{\ell,n}=\frac{\tilde{\omega}_j^{\ell,n}}{\tilde{\omega}_j^{L,n}+\tilde{\omega}_j^{H,n}}, \quad \ell=L,H,
\end{equation}
where
\begin{equation}
\tilde{\omega}_j^{L,n}= \frac{C_L}{\epsilon_t^2}
\end{equation}
is the constant weight associated to the first order non oscillatory solution $u^\IE$ in the time interval $[t^n,t^{n+1}]$, and
\begin{equation}
\tilde{\omega}_j^{H,n}= \frac{C_H}{(\epsilon_t + I_j^3)^2}
\end{equation}
is the weight associated to the \DP\ approximation in the $j$-th cell in the time interval $[t^n,t^{n+1}]$. Here, $\epsilon_t=\Delta t^\tau$, and we always take $\tau=2$ if not otherwise stated.
The $I_j^3$ is a smoothness indicators that measures the regularity of the \DP\ solution. We define $I_j^3$ as contribution of two terms
\begin{equation}\label{e:timeIS3}
I_j^3= I_j^t + I_j^{x,-} + I_j^{x,+},
\end{equation}
where $I_j^t$ and $I_j^{x,\pm}$ are smoothness indicators designed in order to detect discontinuity in time and space, respectively, over the cell $j$. The definition of $I_j^t$ relies on the CE polynomial $P_t^{(3)}$. In fact, at each cell, the CE polynomial changes, and we will have different CE's, and each CE will provide local information on the smoothness of the \DIRK\ advancement in time. We measure the regularity of $P_t^{(3)}$ by the Jiang-Shu smoothness indicator, namely
\begin{equation*}
I_j^t= \sum_{\ell=1}^3 (\Delta t)^{2\ell -1}\int_{t^n}^{t^{n+1}}
\left(\frac{\d^\ell}{\d t^\ell} P^{(3)}_t(t)\right)^2 \d t.
\end{equation*}
Instead, the definition of $I_j^{x,\pm}$ draws inspiration from~\cite{2020Arbogast}, and for the \DIRK[3] method~\eqref{e:dirk3} we have
\begin{align*}
I_j^{x,+} &= (\ca{u}_{j+1}^{\DP,n}-\ca{u}_{j}^{\DP,n})^2 + \sum_{k=1}^2 (\ca{u}_{j+1}^{\DP,(k)}-\ca{u}_{j}^{\DP,(k)})^2 + (\ca{u}_{j+1}^{\DP,n+1}-\ca{u}_{j}^{\DP,n+1})^2 \\
I_j^{x,-} &= (\ca{u}_{j-1}^{\DP,n}-\ca{u}_{j}^{\DP,n})^2 + \sum_{k=1}^2 (\ca{u}_{j-1}^{\DP,(k)}-\ca{u}_{j}^{\DP,(k)})^2 + (\ca{u}_{j-1}^{\DP,n+1}-\ca{u}_{j}^{\DP,n+1})^2, 
\end{align*}
with $\ca{u}_{j}^{\DP,(k)}$ being the approximation at the $k$-th stage.

Finally, we point out that the coefficients $C_L$ and $C_H$ must be carefully chosen. In fact, since the time limited solution~\eqref{e:TimeBlended_general} blends a first order accurate solution with a third order accurate one, according to~\cite{SempliceVisconti:2020} we must choose $C_L=\Delta t^2$ in order to obtain a third order time limited solution.

\subsubsection{Conservative correction}

The two solutions $u^{\IE}$ and $u^{\DP}$ are obtained with conservative schemes, and thus conserve mass. However, the blending~\eqref{e:TimeBlended} itself is not conservative, because it occurs at the cell level, instead of at interfaces.

It is possible therefore that at the $j+\nicefrac12$ interface a mass loss (or gain) is observed. More precisely, the low order predictor can be written as
\begin{equation}
u^{\IE,n+1}_j=\ca{u}_j^n-\frac{\Delta t}{h}\left[ 
\mathcal{F}_{j+1/2}^{\IE}-\mathcal{F}_{j-1/2}^{\IE}
\right],
\end{equation}
while the high order corrector is
\begin{equation}
u^{\DP,n+1}_j=\ca{u}_j^n-\frac{\Delta t}{h}\left[ 
\mathcal{F}_{j+1/2}^{\DP}-\mathcal{F}_{j-1/2}^{\DP}
\right].
\end{equation}
The blended solution therefore is
\begin{equation}
u^{B,n+1}_j=\ca{u}_j^n-\frac{\Delta t}{h}\left[\frac{\omega_j^{H,n}}{C_H}\Delta
\mathcal{F}_{j}^{\DP}+
\left( \omega_j^L- \frac{C_L}{C_H}\omega_j^H\right) \Delta\mathcal{F}_{j}^{\IE}
\right].
\end{equation}
Since the blending is cell-centered, while the fluxes are based on the interfaces, we expect that through the $j+\nicefrac12$ interface there will be a mass loss (or gain) given by
\begin{equation}
\mu_{j+\nicefrac12}^{n+1}= \frac{\Delta t}{C_Hh} \left( (\omega_j^{H,n}-\omega_{j+1}^{H,n})\mathcal{F}_{j+1/2}^{\DP}
+ \left( C_H(\omega_j^{L,n}-\omega_{j+1}^{L,n}) - C_L
(\omega_j^{H,n}-\omega_{j+1}^{H,n})\right) \mathcal{F}_{j+1/2}^{\IE} \right).
\end{equation}
Thus we redistribute the mass lost through the $j+1/2$ interface among the $j$-th and $(j+1)$-th cell, obtaining the limited in time, limited in space Quinpi3 solution, which will be called \QP\ solution. The redistribution is done proportionally to the high order nonlinear weight, so that
\begin{equation}
u_j^{\QP,n+1}=u_j^{B,n+1}+\frac{\omega_j^H}{\omega_j^H+\omega_{j+1}^H} \mu_{j+\nicefrac12} + \frac{\omega_{j}^H}{\omega_j^H+\omega_{j-1}^H} \mu_{j-\nicefrac12}.
\end{equation}
This is the updated solution at time $t^{n+1}$. The resulting scheme is conservative with numerical flux
\begin{align*}
	\mathcal{F}_{j+1/2}^{\QP} =& \frac{1}{C_H}\frac{2\omega_j^{H}\omega_{j+1}^{H}}{\omega_j^{H}+\omega_{j+1}^{H}} \mathcal{F}_{j+1/2}^{\DP} + \\
	&+ \frac{1}{C_H} \left( C_H \frac{\omega_j^L\omega_{j+1}^H+\omega_j^H\omega_{j+1}^L}{\omega_j^H+\omega_{j+1}^H} - C_L\frac{2\omega_j^H\omega_{j+1}^H}{\omega_j^H+\omega_{j+1}^H}\right) \mathcal{F}_{j+1/2}^{\IE} .
\end{align*}
Introducing the "reduced mass" of the high order weights
\[
\tilde{\omega}^H_{j+\nicefrac12} = \frac{2\omega_j^{H}\omega_{j+1}^{H}}{\omega_j^{H}+\omega_{j+1}^{H}},
\]
we can rewrite the conservative flux as
\begin{align*}
	\mathcal{F}_{j+1/2}^{\QP} =& \frac{\tilde{\omega}^H_{j+\nicefrac12}}{C_H} \left[ \mathcal{F}_{j+1/2}^{\DP} + \left( \frac{C_H}2\left( \frac{\omega_j^L}{\omega_j^H}+ \frac{\omega_{j+1}^L}{\omega_{j+1}^H}\right) -C_L\right)  \mathcal{F}_{j+1/2}^{\IE} \right].
\end{align*}
From the formula it is apparent that the high order flux of the interface contributes significantly to the time limited flux only when the cells on both sides are detected as smooth. Moreover, it is also clear that $C_L$ must be infinitesimal.

In~\cite{Marsha:AMR} an analogous conservative correction is employed to ensure the conservation property  at interfaces between grid patches in an Adaptive Mesh Refinement (AMR) algorithm. Other approaches are possible, in particular we refer to~\cite{Ketcheson:fluxbased:2013} which avoids a cell-centered blending by means of flux-based Runge-Kutta. Similar techniques are used in~\cite{2020Arbogast}.

\section{Numerical simulations} \label{sec:numerics}

The purpose of the tests appearing in this section is to study the accuracy of the Quinpi scheme proposed in this work, and to verify their smaller dissipation compared to the first order predictor \IE\ and their improved non oscillatory properties compared to the non-limited in time corrector \DP. Thus we will consider the standard tests which are commonly used in the literature on high order methods for conservation laws: linear advection of non smooth waves, shock formation and interaction in Burgers' equation and the Buckley-Leverett non-convex equation. Furthermore, on one of the tests with singularities, see Figure~\ref{fig:massconservation}, we also demonstrate the need of the conservative correction discussed in the previous section. 

As mentioned in the description of the scheme, when solving nonlinear conservation laws the solution of nonlinear systems is required both for the prediction and the correction step. To this end, we employ the Newton's method. The initial guess to compute the approximation $\overline{u}^{\mathsf{IE},(k)}$ with the predictor at $t^{(k)}=t^n+\theta_k\Delta t$ is chosen as $\overline{u}^{\mathsf{IE},(k-1)}$, $k=1,2,3$. The initial guess to compute the stages $\overline{u}^{(k)}$ of the corrector are the corresponding values $\overline{u}^{\mathsf{IE},(k)}$ of the predictor, for each $k=1,2,3$. The stopping criteria are based on the relative error between two successive approximations and on the norm of the residual. We use a given tolerance $\Delta t^3$, according to the global error of the scheme.

\subsection{Convergence test}

We test the numerical convergence rate of the third-order Quinpi introduced in Section~\ref{sec:scheme} on the nonlinear Burgers' equation
\begin{equation} \label{e:burgers}
	u_t + \left( \frac{u}{2} \right)^2_x = 0,
\end{equation}
with initial condition
\begin{equation} \label{e:smoothIC}
	u_0(x) = 0.5 - 0.25\sin(\pi x)
\end{equation}
on $\Omega = [0,2]$ with periodic domain, and up to the final time $t=1$, i.e.~before the shock appears. The numerical errors, in both $L^1$ and $L^\infty$ norms, and convergence rates are showed in Table~\ref{tab:convergence} for different CFL numbers.

\begin{table}[t!]
	\caption{Orders of convergence of the Quinpi scheme \QP.\label{tab:convergence}}
	\centering
	\subfloat[$\Delta t=h$.]{
		\pgfplotstabletypeset[
		col sep=comma,
		sci zerofill,
		empty cells with={--},
		every head row/.style={before row=\toprule,after row=\midrule},
		every last row/.style={after row=\bottomrule},
		create on use/rate/.style={create col/dyadic refinement rate={1}},
		columns/0/.style={column name={$N$}},
		columns/1/.style={column name={$L^1$ error}},
		columns/rate/.style={fixed zerofill},
		columns={0,1,rate},
		skip rows between index={0}{4}
		]
		{IEDIRK_cou1_err1_constrec.err}
		%
		\pgfplotstabletypeset[
		col sep=comma,
		sci zerofill,
		empty cells with={--},
		every head row/.style={before row=\toprule,after row=\midrule},
		every last row/.style={after row=\bottomrule},
		create on use/rate/.style={create col/dyadic refinement rate={1}},
		columns/1/.style={column name={$L^\infty$ error}},
		columns/rate/.style={fixed zerofill},
		columns={1,rate},
		skip rows between index={0}{4}
		]
		{IEDIRK_cou1_errinf_constrec.err}
	}
	\\
	\subfloat[$\Delta t=10h$.]{
		\pgfplotstabletypeset[
		col sep=comma,
		sci zerofill,
		empty cells with={--},
		every head row/.style={before row=\toprule,after row=\midrule},
		every last row/.style={after row=\bottomrule},
		create on use/rate/.style={create col/dyadic refinement rate={1}},
		columns/0/.style={column name={$N$}},
		columns/1/.style={column name={$L^1$ error}},
		columns/rate/.style={fixed zerofill},
		columns={0,1,rate},
		skip rows between index={0}{3}
		]
		{IEDIRK_cou10_err1_constrec.err}
		%
		\pgfplotstabletypeset[
		col sep=comma,
		sci zerofill,
		empty cells with={--},
		every head row/.style={before row=\toprule,after row=\midrule},
		every last row/.style={after row=\bottomrule},
		create on use/rate/.style={create col/dyadic refinement rate={1}},
		columns/1/.style={column name={$L^\infty$ error}},
		columns/rate/.style={fixed zerofill},
		columns={1,rate},
		skip rows between index={0}{3}
		]
		{IEDIRK_cou10_errinf_constrec.err}
	}
	\\
	\subfloat[$\Delta t=50h$.]{
		\pgfplotstabletypeset[
		col sep=comma,
		sci zerofill,
		empty cells with={--},
		every head row/.style={before row=\toprule,after row=\midrule},
		every last row/.style={after row=\bottomrule},
		create on use/rate/.style={create col/dyadic refinement rate={1}},
		columns/0/.style={column name={$N$}},
		columns/1/.style={column name={$L^1$ error}},
		columns/rate/.style={fixed zerofill},
		columns={0,1,rate},
		skip rows between index={0}{3}
		]
		{IEDIRK_cou50_err1_constrec.err}
		%
		\pgfplotstabletypeset[
		col sep=comma,
		sci zerofill,
		empty cells with={--},
		every head row/.style={before row=\toprule,after row=\midrule},
		every last row/.style={after row=\bottomrule},
		create on use/rate/.style={create col/dyadic refinement rate={1}},
		columns/1/.style={column name={$L^\infty$ error}},
		columns/rate/.style={fixed zerofill},
		columns={1,rate},
		skip rows between index={0}{3}
		]
		{IEDIRK_cou50_errinf_constrec.err}
	}
\end{table}

In the nonlinear blending in time between the low-order solution, i.e.~the composite implicit Euler \IE, and the high-order solution, i.e.~\DP, we use $C_{L} = \Delta t^2$ and $C_{H} = 1-\Delta t^2$ as linear weights. We observe third order  convergence in both norms. In particular, we point out that with large CFL numbers our method reaches smaller errors and faster convergence with respect to the results in~\cite{2020Arbogast} for the same grid.

\begin{figure}[t!]
	\centering
	\includegraphics[width=\textwidth]{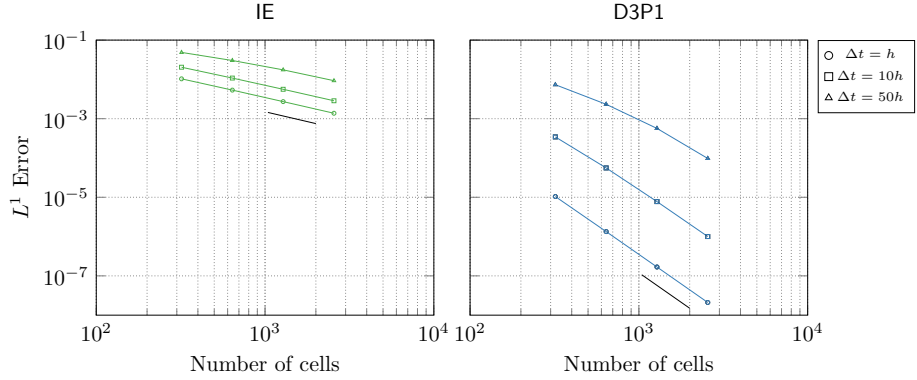}
	\caption{Convergence test of the predictor scheme \IE, left panel, and of the corrector scheme \DP, right panel, for different CFL numbers.\label{fig:convergence}}
\end{figure}

On the same smooth problem, we also test the convergence rate of the low-order predictor \IE\ and the high-order corrector \DP, separately. The convergence tests are shown in Figure~\ref{fig:convergence}. Clearly, with \IE\ we observe first order accuracy. Instead, \DP\ achieves the optimal third order accuracy, as expected.

\subsection{Linear transport}

We consider the linear scalar conservation law
\begin{equation} \label{e:linadv}
u_t+u_x=0
\end{equation}
on the periodic domain in space $\Omega=[-1,1]$, and evolve the initial profile $u_0(x)$ for one period, i.e.~up to final time $t=2$. As initial condition we consider the non-smooth profiles
\begin{subequations}\label{e:nonsmoothIC}
	\begin{equation} \label{e:sindiscontIC}
	u_0(x) = \sin(\pi x) +
	\begin{cases}
	3, & -0.4 \leq x \leq 0.4,\\
	0, & \text{otherwise,}	
	\end{cases}
	\end{equation}
	\begin{equation} \label{e:doublestepIC}
	u_0(x) =
	\begin{cases}
		1, & -0.25 \leq x \leq 0.25,\\
		0, & \text{otherwise.}
	\end{cases}
\end{equation}
\end{subequations}
This problem is used in order to investigate the properties of a scheme to transport non-smooth data with minimal dissipation, dispersion and oscillation effects.

\begin{figure}[t!]
	\centering
	\includegraphics[width=\textwidth]{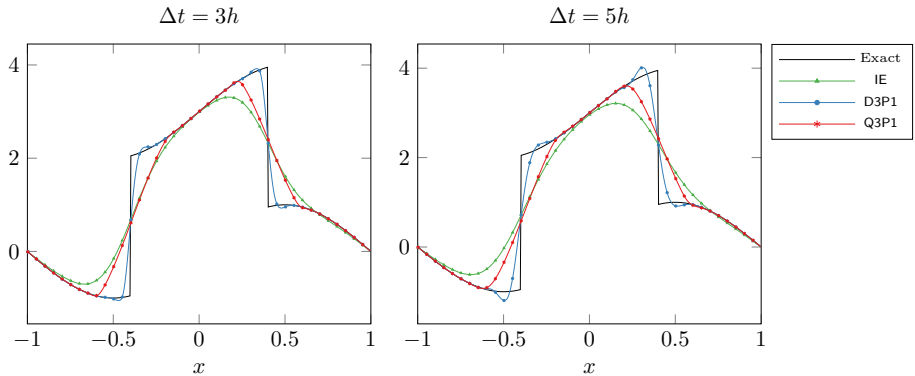}
	\caption{Linear transport equation~\eqref{e:linadv} with initial condition~\eqref{e:sindiscontIC} on $400$ cells at time $t=2$. The markers are used to distinguish the schemes, and are drawn one out of $10$ cells.}
\label{fig:lintra_sindiscont}
\end{figure}

Figure~\ref{fig:lintra_sindiscont} shows the numerical solutions of the linear transport problem with initial profile~\eqref{e:sindiscontIC} computed with the predictor \IE\ and with the corrector methods with no time limiting, i.e.~\DP\, and with blending in time, i.e.~\QP. We consider two different CFL numbers. All the solutions are computed on a grid of $400$ cells. We observe that the low-order predictor is very diffusive, whereas the corrector \DP\ is oscillating across the discontinuities, in particular with CFL number 5. The corrector \QP\, obtained after nonlinear blending of the \IE\ and \DP\ solutions, is much less diffusive than \IE\ and does not produce spurious oscillations, even with large CFL number.

\begin{figure}[t!]
	\centering
	\includegraphics[width=\textwidth]{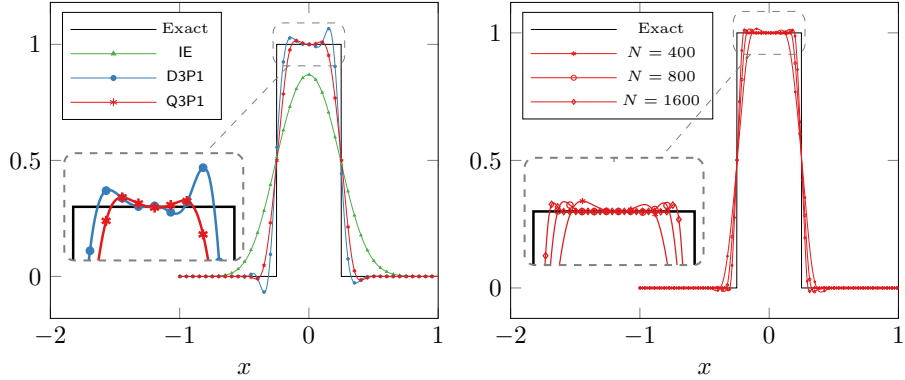}
	\caption{Linear transport equation~\eqref{e:linadv} with initial condition~\eqref{e:doublestepIC} with $\Delta t=5h$ at time $t=2$. The solutions in the left panel are computed on $400$ cells. The right panel shows the solutions obtained with the \QP\ method on different grids. The markers are used to distinguish the schemes, and are drawn one out of $10$ cells on the left and  $15$ cells on the right.\label{fig:lintra_doublestep}}
\end{figure}

In Figure~\ref{fig:lintra_doublestep} we provide the numerical solutions of the linear transport problem with the initial double-step profile~\eqref{e:doublestepIC} with a zoom on the top part of the non-smooth region of the solution. The simulations are performed with $\Delta t= 5h$. In the left panel, we compare the three methods on a grid of $400$ cells. Also in this test we observe that the novel method \QP\ introduced in this work presents less oscillations than \DP\ close to discontinuities, and it is less dissipative than \IE. The right panel shows the approximation provided by \QP\ on different grids. As we expect, on finer grids the frequency of the oscillations increases, whereas the amplitude slightly decreases.

\subsection{Burgers' equation}

We investigate the behavior of the schemes on the nonlinear Burgers' equation~\eqref{e:burgers} for different initial conditions.

\subsubsection{Smooth profile: Shock formation}

As in~\cite{2020Arbogast}, we consider the smooth initial condition~\eqref{e:smoothIC} on the periodic domain $\Omega=[0,2]$, and up to the final time $t=2$, i.e.~after shock formation.

\begin{figure}[t!]
	\centering
	\includegraphics[width=\textwidth]{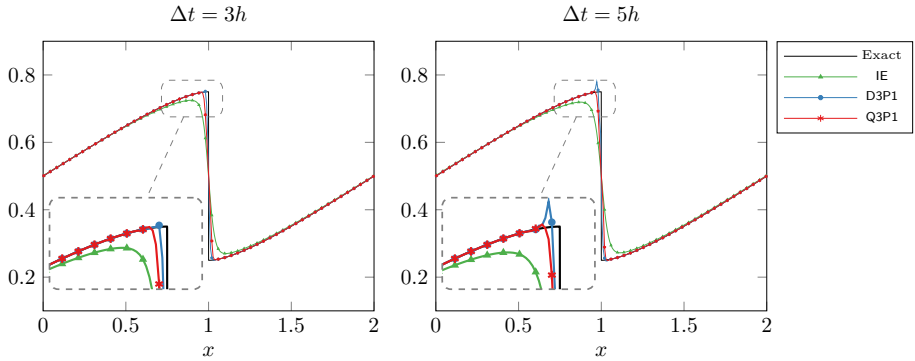}
	\caption{Burgers' equation~\eqref{e:burgers} with initial condition~\eqref{e:smoothIC} on $N=256$ cells at time $t=2$. The markers are used to distinguish the schemes, and are drawn each $5$ cells.\label{fig:burgers_arbogast}}
\end{figure}

The results in Fig.~\ref{fig:burgers_arbogast} are obtained with $256$ cell with two values of the CFL number. In both cases the \QP\ scheme is slightly more diffusive than its non-blended version \DP, and exhibits a much lower dissipation than the first order predictor \IE. For $\Delta t = 3h$, all schemes do not produce large oscillations near the discontinuity, and the solutions of \QP\ and \DP\ are very close. We appreciate the difference when $\Delta t = 5h$. In fact, the \QP\ scheme reduces the oscillations created by \DP, while maintaining a very high resolution.

\begin{figure}[t!]
	\centering
	\includegraphics[width=\textwidth]{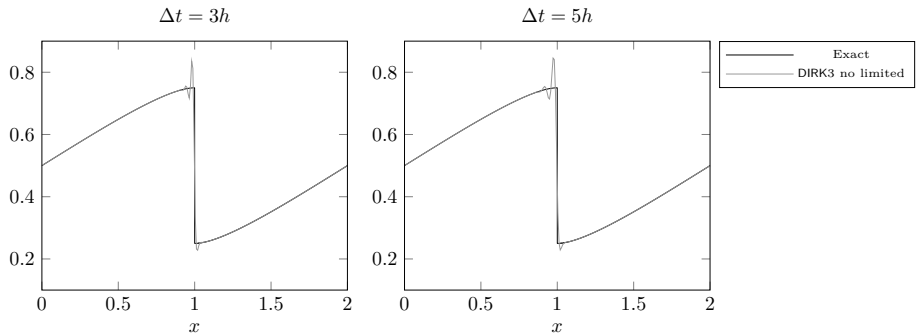}
	\caption{Numerical solution of the Burgers' equation~\eqref{e:burgers} with initial condition~\eqref{e:smoothIC} on $N=256$ cells at time $t=2$ obtained with the space-time non-limited scheme.\label{fig:burgers_arbogast_noLimiting}}
\end{figure}

On this test, we show also the numerical approximation provided by the space-time non-limited scheme, cf.~Fig.~\ref{fig:burgers_arbogast_noLimiting}. The setup of the simulation is as in Fig.~\ref{fig:burgers_arbogast}, i.e.~we consider $256$ cells with CFL numbers $3$ and $5$. Compared to Fig.~\ref{fig:burgers_arbogast}, we observe the importance of the limiting technique in order to avoid the very large spurious oscillations appearing also at moderate CFL numbers.

\subsubsection{Non-smooth profile}

We test the Burgers' equation on the discontinuous initial condition~\eqref{e:doublestepIC}, on the periodic domain $\Omega=[-1,1]$, and up to the final time $t=0.5$.

\begin{figure}[t!]
	\centering
	\includegraphics[width=\textwidth]{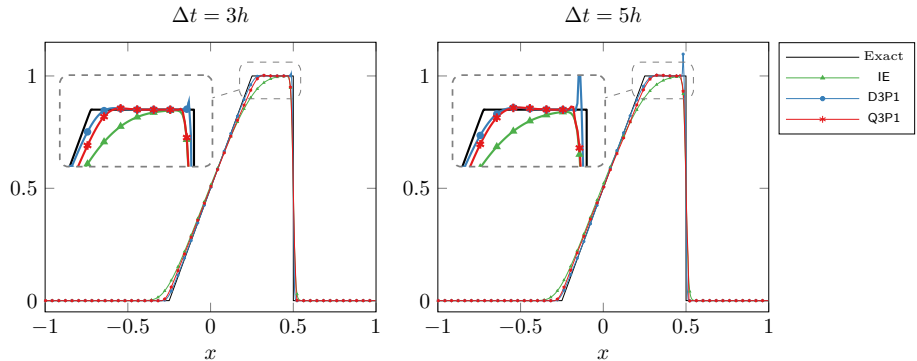}
	\caption{Burgers' equation~\eqref{e:burgers} with initial condition~\eqref{e:doublestepIC} on $400$ cells at time $t=0.5$. The markers are used to distinguish the schemes, and are drawn every $8$ cells.\label{fig:burgers_doublestep}}
\end{figure}

The numerical solutions are shown in Figure~\ref{fig:burgers_doublestep} with $400$ cells, with $\Delta t=3h$ and $\Delta t=5h$. The nonlinear Burgers' equation develops a rarefaction and a moving shock. Again, we observe the ability of the new implicit scheme \QP\ of increasing the accuracy on smooth zones compared to \IE, and, at the same time, reducing the spurious oscillations across the shock.

\begin{figure}[t!]
	\centering
	\includegraphics[width=\textwidth]{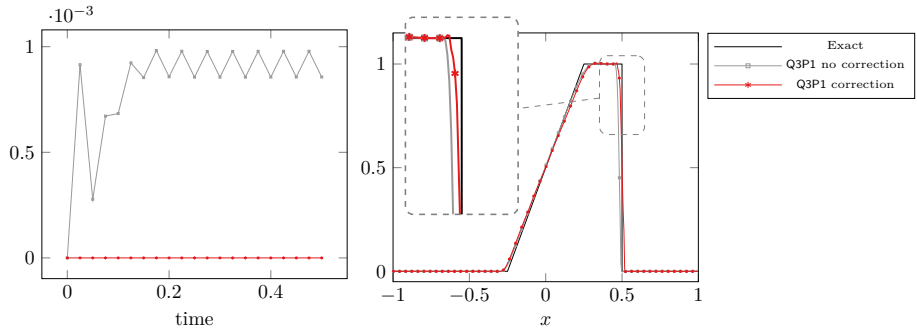}
	\caption{Right plot: numerical solutions of \QP, with and without conservative correction, on the Burgers' equation~\eqref{e:burgers} with initial condition~\eqref{e:doublestepIC} using $400$ cells at time $t=5$ and with $\Delta t=5h$. The markers are used to distinguish the schemes, and are drawn every $8$ cells. Left plot: deviation in time of the total mass of the numerical solution from the exact total mass.\label{fig:massconservation}}
\end{figure} 

On this particular test, we provide a numerical evidence of the need and of the effectiveness of the conservative correction introduced after the nonlinear blending in time between the predictor scheme \IE\ and the third-order corrector \DP. In the right panel of Figure~\ref{fig:massconservation} we show the solutions provided by the Quinpi scheme \QP\ with and without the conservative correction, using $400$ cells and $\Delta t=5h$. Instead, in the left panel of Figure~\ref{fig:massconservation} we show the behavior in time of the deviation of the mass of the numerical solution from the mass of the initial condition. We observe that, without correction, the \QP\ scheme does not capture the correct shock location because of the mass lost. The conservative correction allows to predict the shock at the correct location and the mass is conserved at all times. 

\begin{figure}[t!]
	\centering
	\includegraphics[width=\textwidth]{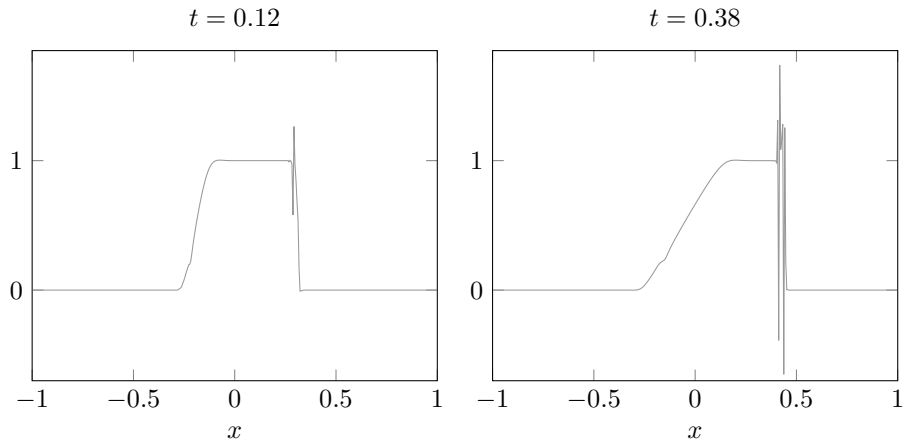}
	\caption{Numerical solution of the Burgers' equation~\eqref{e:burgers} with initial condition~\eqref{e:doublestepIC} on $400$ cells with $\Delta t=5h$ obtained with an explicit first order predictor.\label{fig:burgers_doublestep_expPred}}
\end{figure}

Finally, in Figure~\ref{fig:burgers_doublestep_expPred} we show the numerical approximations obtained with a scheme using an explicit first order predictor, on a grid of $400$ cells with $\Delta t=5h$. The solutions are provided at two different times, and we observe that the use of an explicit predictor is not enough to prevent spurious oscillations in the high order scheme, at relatively high CFL numbers. Contrary to~\cite{Gottlieb:iWENO:2006}, in this numerical test we compute the explicit predictor at each intermediate time defined by the abscissae of the \DIRK\ scheme.

\subsubsection{Shock interaction}

We consider Burgers' equation with smooth initial condition
\begin{equation} \label{e:shockinteractionIC}
u_0(x) = 0.2 -\sin(\pi x) + \sin(2 \pi x) 
\end{equation}
on the periodic domain $\Omega=[-1,1]$, and with $\Delta t = 5h$. This test allows to compare the behavior of the schemes on both shock formation and shock interaction. In fact, the exact solution is characterized by the formation of two shocks which collide at a larger time, developing a single discontinuity.

\begin{figure}[t!]
	\centering
	\includegraphics[width=\textwidth]{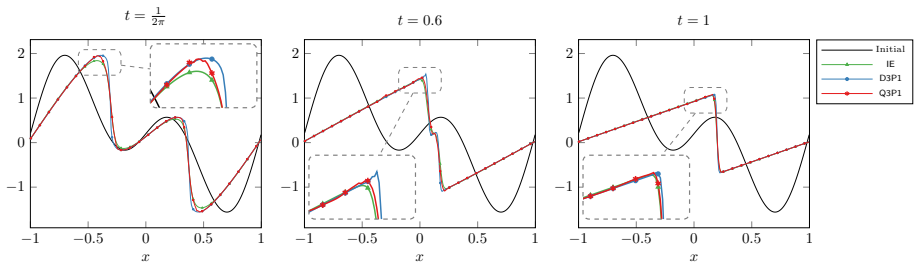}
\caption{Burgers' equation~\eqref{e:burgers} with initial condition~\eqref{e:shockinteractionIC} on $256$ cells with $\Delta t=5h$, at three different times. The markers are used to distinguish the schemes, and are drawn each $10$ cells.}
\label{fig:burgers_shockinteraction}
\end{figure}

In Figure~\ref{fig:burgers_shockinteraction} we show the numerical solutions at three snapshots: at $t=\frac{1}{2\pi}$, when the two shocks occur, at $t=0.6$, which is slightly before the interaction of the two shocks, and finally at $t=1$, shortly after the shock interaction. It is clear that \QP\ does not produce spurious oscillations, and its profile has a higher resolution with respect to \IE.

\subsection{Buckley-Leverett equation}

The last numerical tests are performed on the Buckley-Leverett equation
\begin{equation} \label{e:buckleyleverett}
	u_t + \left( \frac{u^2}{u^2+\frac13(1-u)^2} \right)_x = 0
\end{equation}
which is characterized by a non-convex flux function. We consider the same setup as in~\cite{2020Arbogast}. Therefore, the initial profile is the step function
\begin{equation} \label{e:stepIC}
u_0(x) =
\begin{cases}
0.5, & -0.25 \leq x \leq 0.25,\\
0, & \text{otherwise.}
\end{cases}
\end{equation}
on the periodic domain $\Omega=[0,1]$, and up to the final time $t=0.085$.

\begin{figure}[t!]
	\centering
	\includegraphics[width=\textwidth]{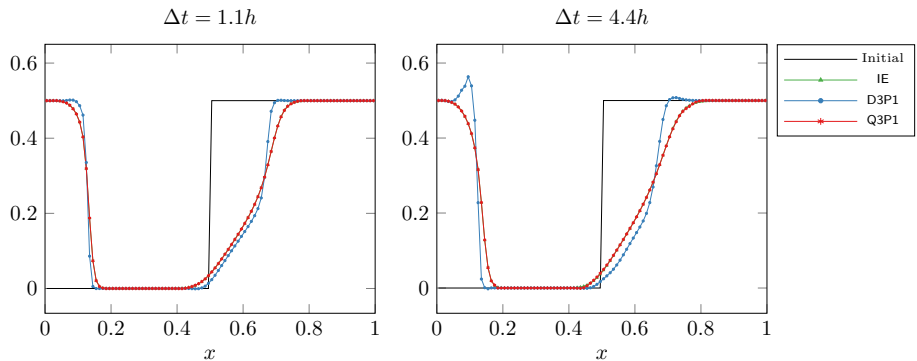}
	\caption{Buckley-Leverett equation~\eqref{e:buckleyleverett} 
	with initial condition~\eqref{e:stepIC} on $100$ cells at 
	time $t=0.085$.}
\label{fig:buckleyleverett}
\end{figure}

In Figure~\ref{fig:buckleyleverett} we show the results on a grid of $100$ cells with $\Delta t = 1.1h$ and $\Delta t = 4.4h$. For this particular example we set $\epsilon_t=\Delta t^3$ in the nonlinear blending in time. The results produced by the three schemes are very similar when the CFL is small. With large CFL, instead, we observe that \QP\ improves \DP\ avoiding the spurious overshoots.

\subsection{Performance of the third order Quinpi scheme}

The implicit method proposed in this paper is clearly more expensive comparing to explicit schemes. On the contrary, the benefit of an implicit scheme is to achieve larger CFL numbers and more efficient computations for stiff problems. It is thus important to identify the regimes where it is convenient using an implicit method. To this end, in the following we compare the computational CPU times required by each time step of the explicit third order SSP Runge-Kutta scheme and of the \QP\ scheme.
In Table~\ref{tab:cpu} we report the results obtained both on the linear equation~\eqref{e:linadv} on $[-1,1]$
and on the Burgers' equation~\eqref{e:burgers} on $[0,2]$ with initial condition~\eqref{e:shockinteractionIC}.
To make the test reliable and fair, the parameters are chosen to have a comparable and large enough total time of execution of the methods. The results are obtained on a
quadcore Intel Core i7-6600U with clock speed 2.60GHz.
\begin{table}[t!]
	\caption{CPU times for each step of the SSP-RK3 explicit scheme and of the \QP\ implicit scheme on linear and nonlinear problems.\label{tab:cpu}}
	\centering
	\subfloat[Linear problem.]{
		\begin{tabular}{c|c|c|c}
			\hline
			Cells $N$ & SSP-RK3 & Q3P1 & Ratio $r_N$ \\
			\hline
			$200$ & 0.0021 & 0.0062 & 2.95 \\
			$400$ & 0.0033 & 0.0093 & 2.82 \\
			$800$ & 0.0061 & 0.0155 & 2.54 \\
			$1600$ & 0.095 & 0.0293 & 3.08 \\
			\hline
		\end{tabular}
	}
	\\
	\subfloat[Nonlinear problem before shock formation.]{
		\begin{tabular}{c|c|c|c}
			\hline
			Cells $N$ & SSP-RK3 & Q3P1 & Ratio $r_N$\\
			\hline
			$200$ & 0.0023 & 0.0103 & 4.48 \\
			$400$ & 0.0032 & 0.0148 & 4.62 \\
			$800$ & 0.0066 & 0.0229 & 3.47 \\
			$1600$ & 0.0091 & 0.0361 & 3.97 \\
			\hline
		\end{tabular}
	}
	\subfloat[Nonlinear problem after shock formation.]{
		\begin{tabular}{c|c|c|c}
			\hline
			Cells $N$ & SSP-RK3 & Q3P1 & Ratio $r_N$ \\
			\hline
			$200$ & 0.0023 & 0.0107 & 4.65 \\
			$400$ & 0.0033 & 0.0158 & 4.79 \\
			$800$ & 0.0068 & 0.0255 & 3.75 \\
			$1600$ & 0.0095 & 0.0566 & 5.96 \\
			\hline
		\end{tabular}
	}
\end{table}
As it is expected, a step of the implicit scheme is more expensive than a step of the explicit one and the former is more convenient when the problem does not require a time step $\Delta t \leq C \frac{h}{\lambda_{\max} \left\lceil r_N \right\rceil}$, where $\lambda_{\max} = \max_u |f'(u)|$.

\begin{figure}[t!]
	\centering
	\subfloat[$400$ cells.]{
		\includegraphics[width=\textwidth]{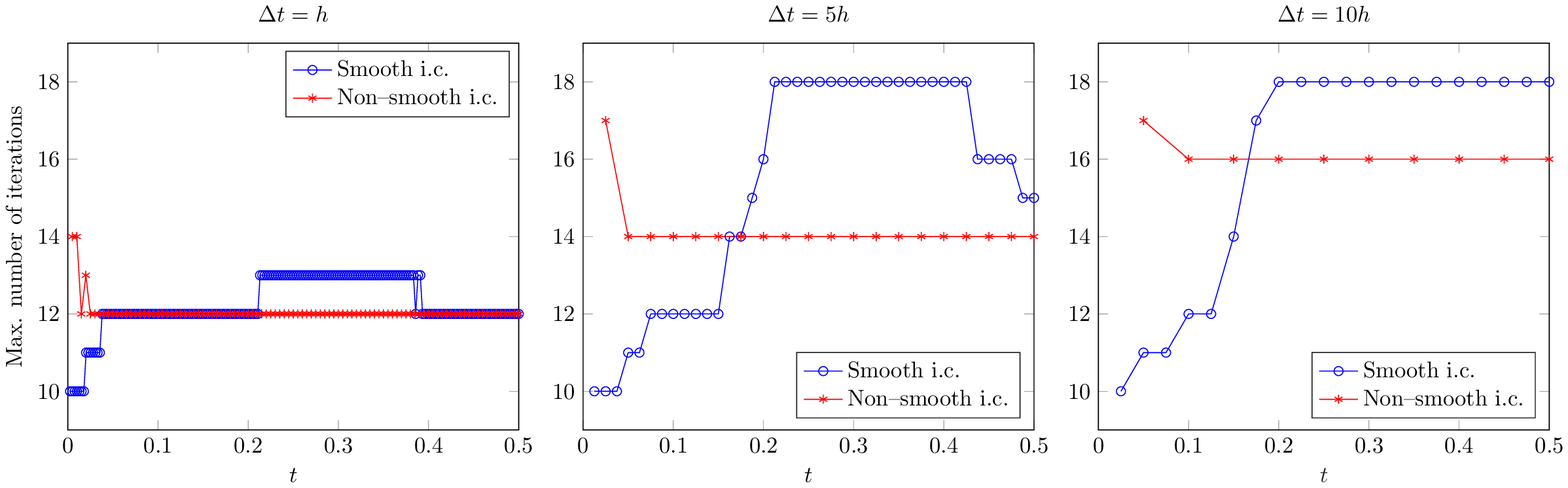}
	}
	\\
	\subfloat[$800$ cells.]{
		\includegraphics[width=\textwidth]{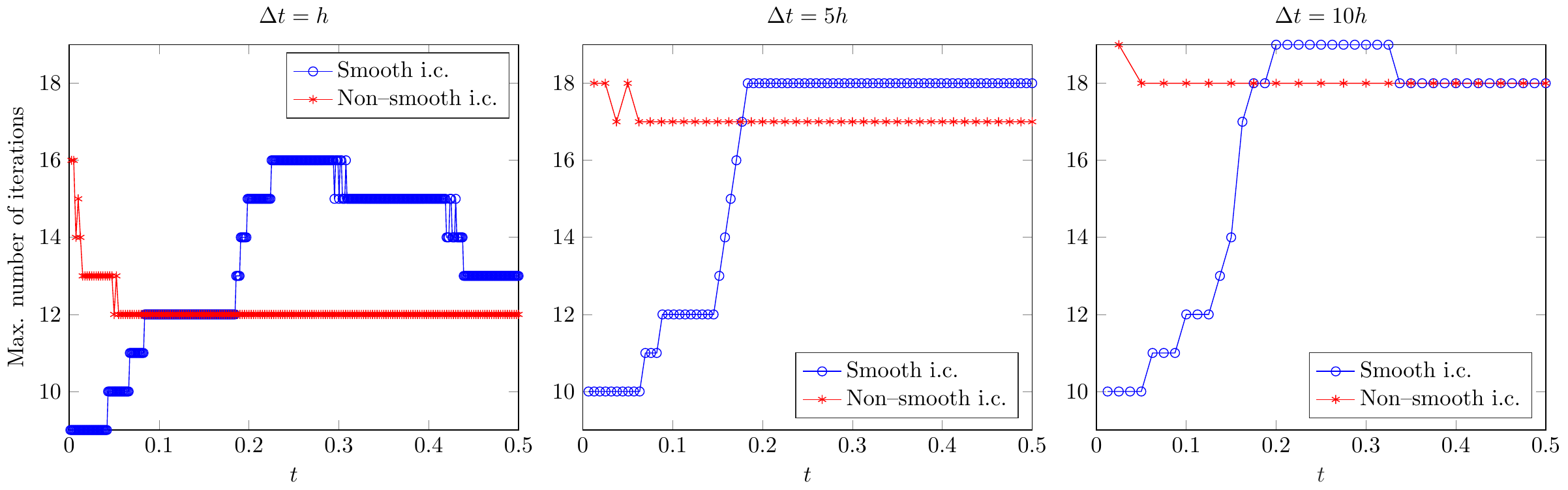}
	}
	\caption{Total number of iterations in each time step required for the convergence of the Newton's method when solving the Burgers' equation up to $t=0.5$, with initial condition~\eqref{e:shockinteractionIC} (blue lines) and~\eqref{e:doublestepIC} (red lines). Three different values of the Courant number are considered.\label{fig:newton_convergence}}
\end{figure}

The complexity of the implicit scheme increases particularly on the nonlinear problem. In fact, the \QP\ scheme requires the solution of the Newton's method six times in each time step. However, the number of iterations remains limited and the convergence of the method is not problematic even in presence of shocks. In Figure~\ref{fig:newton_convergence} we show the total number of iterations required for the Newton's method in each time step in solving the Burgers' equation~\eqref{e:burgers} with smooth~\eqref{e:shockinteractionIC} and non--smooth~\eqref{e:doublestepIC} initial profiles	up to $t=0.5$ and three different Courant numbers. The top panels are obtained with $400$ cells, whereas the bottom panels are with $800$ cells.
For the smooth initial condition, we observe that the number of iterations is not larger than 2 for each Newton's method up to a time when the solution remains smooth. Instead, the number of iterations slightly increases when the solution becomes discontinuous. In fact, each Newton's method is converging with maximum 3 iterations. For the double step initial condition, the number of iterations is large, but not more than 3 for each Newton's method, in the first time step and decreases when the rarefaction appears.

\section{Conclusions}\label{sect:conclusion}

In this work, we have proposed a new approach to the integration of hyperbolic conservation laws with high order implicit schemes. The main characteristic of this framework is to use low order implicit predictors with a double purpose. First, the predictor is used to determine the nonlinear weights in a \CWENO\ or \WENO\ high order space reconstruction. In this fashion, one greatly simplifies the differentiatiation of  the \WENO\ weights when computing the Jacobian of the numerical fluxes. Second, the predictor is used as low order approximation of the solution and is blended with the high order solution in order to achieve limiting also in time.

The resulting scheme is {\em linear} with respect to the solution at the new time level on linear equations, unlike most, if not all,  high order implicit schemes available. The non-linearity of the scheme is linked only to the non-linearity of the flux. This does not mean that the coefficients appearing in the scheme are constant. It means  that the non-linearities in the limiting in space and time of the scheme involve only the predictor, which is already known when the high order solution is evolved in time.

We expect the new scheme to have applications in many stiff problems, as low Mach gas dynamics or kinetic problems. Future work will involve the application of the Quinpi approach to stiff gas dynamics, the exploration of this new framework with BDF time integration, and extensions to higher order.

\paragraph{Acknowledgments} This work was partly supported by MIUR (Ministry of University and Research) PRIN2017 project number 2017KKJP4X.

No conflict of interest is extant in the present work.

\bibliographystyle{plain}      
\bibliography{Quinpi}

\begin{thebibliography}{10}

\bibitem{2017AbbateAllSpeed}
E.~Abbate, A.~Iollo, and G.~Puppo.
\newblock An all-speed relaxation scheme for gases and compressible materials.
\newblock {\em J. Comput. Phys.}, 351:1--24, 2017.

\bibitem{Alexander1977}
R.~Alexander.
\newblock Diagonally implicit {R}unge-{K}utta methods for stiff {O.D.E.}'s.
\newblock {\em {SIAM} J. Numer. Anal.}, 14(6):1006--1021, 1977.

\bibitem{2020Arbogast}
T.~Arbogast, C.S. Huang, X.~Zhao, and D.~N. King.
\newblock A third order, implicit, finite volume, adaptive {R}unge-{K}utta
  {WENO} scheme for advection-diffusion equations.
\newblock {\em Comput. Methods Appl. Mech. Engrg.}, 368, 2020.

\bibitem{Balsara:AOWENO}
D.~S. Balsara, S.~Garain, and C.~W. Shu.
\newblock An efficient class of {WENO} schemes with adaptive order.
\newblock {\em J. Comput. Phys.}, 326:780--804, 2016.

\bibitem{Marsha:AMR}
M.~J. Berger and R.~J. Le~Veque.
\newblock Adaptive mesh refinement using wave-propagation algorithms for
  hyperbolic systems.
\newblock {\em {SIAM} J. Numer. Anal.}, 35(6):2298--2316, 1998.

\bibitem{2018BoscarinoScandurra}
S.~Boscarino, G.~Russo, and L.~Scandurra.
\newblock All {M}ach number second order semi-implicit scheme for the {E}uler
  equations of gas dynamics.
\newblock {\em Journal of Scientific Computing}, 77(2):850--884, 2018.

\bibitem{CWENOandaluz}
M.~J. Castro-D{\`i}az and M.~Semplice.
\newblock Third- and fourth-order well-balanced schemes for the shallow water
  equations based on the {CWENO} reconstruction.
\newblock {\em Int. J. Numer. Meth. Fluid}, 89(8):304--325, 2019.

\bibitem{CPSV:cweno}
I.~Cravero, G.~Puppo, M.~Semplice, and G.~Visconti.
\newblock {CWENO}: uniformly accurate reconstructions for balance laws.
\newblock {\em Math. Comp.}, 87(312):1689--1719, 2018.

\bibitem{CS:epsweno}
I.~Cravero and M.~Semplice.
\newblock On the accuracy of {WENO} and {CWENO} reconstructions of third order
  on nonuniform meshes.
\newblock {\em Journal of Scientific Computing}, 67:1219--1246, 2016.

\bibitem{CSV19:cwenoz}
I.~Cravero, M.~Semplice, and G.~Visconti.
\newblock Optimal definition of the nonlinear weights in multidimensional
  {C}entral {WENOZ} reconstructions.
\newblock {\em {SIAM} J. Numer. Anal.}, 57(5):2328--2358, 2019.

\bibitem{2011DegondTang}
P.~Degond and M.~Tang.
\newblock All speed scheme for the low {M}ach number limit of the isentropic
  {E}uler equations.
\newblock {\em Comm. Computat. Phys.}, 10(1):1--31, 2011.

\bibitem{2010DellacherieLowMach}
S.~Dellacherie.
\newblock Analysis of {G}odunov type schemes applied to the compressible
  {E}uler system at low {M}ach number.
\newblock {\em J. Comput. Phys.}, 229(4):978--1016, 2010.

\bibitem{2017DimarcoLoubere_LowMachAP}
G.~Dimarco, R.~Loubere, and M.~H. Vignal.
\newblock Study of a new asymptotic preserving scheme for the {E}uler system in
  the low {M}ach number limit.
\newblock {\em {SIAM} J. Sci. Comput.}, 39(5):A2099--A2128, 2017.

\bibitem{DimarcoPareschiActa}
G.~Dimarco and L.~Pareschi.
\newblock Numerical methods for kinetic equations.
\newblock {\em Acta Numerica}, 23:369--520, 2014.

\bibitem{Dumbser:2008:PnPm}
M.~Dumbser, D.~S. Balsara, E.~F. Toro, and C.-D. Munz.
\newblock A unified framework for the construction of one-step finite volume
  and discontinuous {G}alerkin schemes on unstructured meshes.
\newblock {\em J. Comput. Phys.}, 227:8209--8253, 2008.

\bibitem{DBSR:ADER_CWENO}
M.~Dumbser, W.~Boscheri, M.~Semplice, and G.~Russo.
\newblock Central weighted {ENO} schemes for hyperbolic conservation laws on
  fixed and moving unstructered meshes.
\newblock {\em {SIAM} J. Sci. Comput.}, 39(6):A2564--A2591, 2017.

\bibitem{2007DurasaisamyBaeder}
K.~Duraisamy and J.D. D.~Baeder.
\newblock Implicit scheme for hyperbolic conservation laws using non
  oscillatory reconstruction in space and time.
\newblock {\em {SIAM} J. Sci. Comput.}, 29:2607--2620, 2007.

\bibitem{2003DurasaisamyBaeder}
K.~Duraisamy, J.D. D.~Baeder, and J.~G. Liu.
\newblock Concepts and application of time-limiters to high resolution schemes.
\newblock {\em Journal of Scientific Computing}, 19:139--162, 2003.

\bibitem{2001Forth}
S.~A. Forth.
\newblock A second order accurate, space-time limited, {BDF} scheme for the
  linear advection equation.
\newblock In {\em {G}odunov Methods}, pages 335--342. Springer, 2001.

\bibitem{Gottlieb:iWENO:2006}
S.~Gottlieb, J.~S. Mullen, and S.~J. Ruuth.
\newblock {A Fifth Order Flux Implicit WENO Method}.
\newblock {\em Journal of Scientific Computing}, 27:271--287, 2006.

\bibitem{2001Gottlieb}
S.~Gottlieb, C.W. Shu, and E.~Tadmor.
\newblock Strong stability preserving high-order time discretization methods.
\newblock {\em SIAM Rev.}, 43:73--85, 2001.

\bibitem{1983Harten_HighResolution}
A.~Harten.
\newblock High resolution schemes for hyperbolic conservation laws.
\newblock {\em J. Comput. Phys.}, 49(3):357--393, 1983.

\bibitem{1984Harten}
A.~Harten.
\newblock On a class of high resolution total-variation-stable
  finite-difference schemes.
\newblock {\em {SIAM} J. Numer. Anal.}, 21, 1984.

\bibitem{JiangShu:96}
G.-S. Jiang and C.-W. Shu.
\newblock Efficient implementation of weighted {ENO} schemes.
\newblock {\em J. Comput. Phys.}, 126:202--228, 1996.

\bibitem{1995JinXin_Relaxation}
S.~Jin and Z.~Xin.
\newblock The relaxation schemes for systems of conservation laws in arbitrary
  space dimensions.
\newblock {\em Communications on Pure and Applied Mathematics}, 48(3):235--276,
  1995.

\bibitem{Ketcheson:fluxbased:2013}
D.~I. Ketcheson, C.~B. MacDonald, and S.~J. Ruuth.
\newblock Spatially partitioned embedded {R}unge-{K}utta methods.
\newblock {\em {SIAM} J. Numer. Anal.}, 51(5):2887--2910, 2013.

\bibitem{Kolb:14}
O.~Kolb.
\newblock On the full and global accuracy of a compact third order {WENO}
  scheme.
\newblock {\em {SIAM} J. Numer. Anal.}, 52(5):2335--2355, 2014.

\bibitem{LeVeque:book}
R.~Le~Veque.
\newblock {\em Finite Volume Methods for Hyperbolic Problems}.
\newblock Cambridge Texts in Applied Mathematics. Cambridge University Press,
  2004.

\bibitem{MicroMacro}
M.~Lemou and L.~Mieussens.
\newblock A new asymptotic preserving scheme based on micro–macro formulation
  for linear kinetic equations in the diffusion limit.
\newblock {\em SIAM J. Sci. Comput.}, 31:334--368, 2008.

\bibitem{LPR:00:SIAMJSciComp}
D.~Levy, G.~Puppo, and G.~Russo.
\newblock Compact central {WENO} schemes for multidimensional conservation
  laws.
\newblock {\em {SIAM} J. Sci. Comput.}, 22(2):656--672, 2000.

\bibitem{NorsettWanner:1981}
S.~P. N{\o}rsett and G.~Wanner.
\newblock Perturbed collocation and {R}unge-{K}utta methods.
\newblock {\em Numer. Math.}, 38:193--208, 1981.

\bibitem{MiMe}
S.~Pieraccini and G.~Puppo.
\newblock Microscopically implicit–macroscopically explicit schemes for the
  {BGK} equation.
\newblock {\em J. Comput. Phys.}, 231:299--327, 2012.

\bibitem{SCR:CWENOquadtree}
M.~Semplice, A.~Coco, and G.~Russo.
\newblock Adaptive mesh refinement for hyperbolic systems based on third-order
  compact {WENO} reconstruction.
\newblock {\em Journal of Scientific Computing}, 66:692--724, 2016.

\bibitem{SempliceVisconti:2020}
M.~Semplice and G.~Visconti.
\newblock Efficient implementation of adaptive order reconstructions.
\newblock {\em Journal of Scientific Computing}, 83:6, 2020.

\bibitem{Shu97}
C.~W. Shu.
\newblock {Essentially non-oscillatory and weighted essentially non-oscillatory
  schemes for hyperbolic conservation laws}.
\newblock In {\em {Advanced numerical approximation of nonlinear hyperbolic
  equations (Cetraro, 1997)}}, volume 1697 of {\em {Lecture Notes in Math.}},
  pages 325--432. Springer, Berlin, 1998.

\bibitem{1984Sweby_TVD}
P.K. Sweby.
\newblock High resolution schemes using flux limiters for hyperbolic
  conservation laws.
\newblock {\em {SIAM} J. Numer. Anal.}, 21(5):995--1011, 1984.

\bibitem{2017Tavelli_SemiImplicitAllMach}
M.~Tavelli and M.~Dumbser.
\newblock A pressure-based semi-implicit space-time discontinuous {G}alerkin
  method on staggered unstructured meshes for the solution of the compressible
  {N}avier-{S}tokes equations at all {M}ach numbers.
\newblock {\em Journal of Computational Physics}, 341:341--376, 2017.

\bibitem{Zennaro}
M.~Zennaro.
\newblock {Natural Continuous Extensions of Runge-Kutta Methods}.
\newblock {\em Math. Comp.}, 46:119--133, 1986.

\end{thebibliography}

\end{document}